\documentclass[12pt,a4paper]{article}
\usepackage{amsmath,amsfonts,amssymb,amscd}
\usepackage{xypic}

\usepackage[unicode=true,
bookmarks=true,bookmarksnumbered=true,bookmarksopen=false,
breaklinks=false,pdfborder={0 0 1},backref=false,colorlinks=true]
{hyperref}

\usepackage[a4paper]{geometry}
\def\Ker{\operatorname{Ker}}

\def\id{\operatorname{id}}

\def\Im{\operatorname{Im}}

\newcounter{th}
\def\t{\refstepcounter{th}{\bf \noindent{Theorem} \arabic{th}. }}

\newcounter{prop}
\def\prop{\refstepcounter{prop}{\bf \noindent{Proposition} \arabic{prop}. }}

\newcounter{lem}

\newcounter{de}

\newcounter{ex}

\begin{document}

\begin{center}
    {\LARGE{\bf A geometric approach to  $1$-singular Gelfand-Tsetlin $\mathfrak {gl}_n$-modules}}
    
\end{center}
\medskip

\medskip

\begin{center}
	{\Large  Elizaveta Vishnyakova}
\end{center}

\bigskip

\begin{abstract}	
This paper is devoted to an elementary new construction of $1$-singular  Gelfand-Tsetlin modules using complex geometry. We introduce a universal ring $\mathcal D_v$ together with the vector space $\mathcal S=\mathcal S(\mathcal D_v)$ with  basis $\mathcal B_v = \mathcal B(\mathcal D_v)$ consisted of some local distributions such that $\mathcal S$ is a natural $\mathcal D_v$-module. For any homomorphism of rings $\mathcal U(\mathfrak{h}) \to \mathcal D_v$, where $\mathfrak{h}$ is a Lie algebra, it follows that  $\mathcal S$ is also an $\mathfrak{h}$-module. We observe that the homomorphism of rings constructed in \cite{FO} is a homomorphism of type $\mathcal U(\mathfrak{gl}_n(\mathbb C)) \to \mathcal D_v$. Using this observation we obtain a construction of the universal $1$-singular Gelfand-Tsetlin $\mathfrak{gl}_n(\mathbb C)$-module from \cite{Futorny}.

\end{abstract}

\section{Introduction}

This paper is devoted to a new elementary geometric construction of the universal $1$-singular Gelfand-Tsetlin module. Denote $\mathfrak g_k:=\mathfrak{gl}_k(\mathbb C)$, where $k=1,\ldots n$, and consider the flag $\mathfrak g_1\subset \mathfrak g_2 \subset\cdots \subset \mathfrak g_{n-1} \subset \mathfrak g_n$ of Lie algebras, where $\mathfrak g_{k-1} \subset \mathfrak g_k$ is the inclusion with respect to the left top corner. This flag gives rise to the following flag of universal enveloping algebras 
$$
\mathcal U(\mathfrak g_1)\subset \mathcal U(\mathfrak g_2) \subset\cdots \subset \mathcal U(\mathfrak g_{n-1}) \subset \mathcal U(\mathfrak g_n). 
$$
Denote by $\mathcal Z_k$ the center of $\mathcal U(\mathfrak g_{k})$. Then
the subalgebra $\Gamma\subset \mathcal U(\mathfrak g_n)$, generated by $\mathcal Z_k$, where $k = 1,  \ldots, n$, is a maximal commutative subalgebra \cite{Ov1}. It is called the {\it Gelfand-Tsetlin subalgebra}. A $\mathcal U(\mathfrak g_n)$-module $M$ is called a  {\it Gelfand-Tsetlin  module} if the action of $\Gamma$ on $M$ is locally finite.

In the classical Gelfand-Tsetlin theory \cite{Gelfand} an explicit construction of an action of $\mathfrak g_n$ with respect to a basis consisting of Gelfand-Tsetlin tableaux is given providing explicit formulas for $\mathfrak g_n$-action. These formulas for $\mathfrak g_n$-action are called {\it classical Gelfand-Tsetlin formulas}. It was noticed in \cite{DFO1,DFO2,DFO3,DFO4} that the classical Gelfand-Tsetlin formulas may be used to obtain a family of infinite dimensional Gelfand-Tsetlin  modules: so-called {\it generic regular Gelfand-Tsetlin modules}. An essential progress in the theory of Gelfand-Tsetlin  modules was done in \cite{Ov1,Ov2} and later in  \cite{FO}. In particular the following important construction was obtained there. Let $V\simeq \mathbb C^{n(n+1)/2}$ be the vector space of all Gelfand-Tsetlin tableaux of fixed order $n$, see the main text for details. Denote by $\gimel$ a certain abelian group acting freely on $V$ and by $\mathcal M\star \gimel$ the sheaf of meromorphic functions on $V$ with values in $\gimel$. Then there exists a ring structure on $\mathcal R:=H^0(V,\mathcal M \star \gimel)$  such that 
 the classical Gelfand-Tsetlin formulas define a ring homomorphism $\Phi:	\mathcal U(\mathfrak g_n) \to \mathcal R$.  In the case when $\Im\Phi$  is holomorphic at a neighborhood of the orbit $\gimel(v)$ of a point $v\in V$, we can define a $\mathfrak g_n$-module structure on the vector space with the basis $\{ev_v, \,\, v\in  \gimel(v)\}$, where $ev_v$ is the evaluation map at the point $v$. These $\mathfrak g_n$-modules are exactly generic regular Gelfand-Tsetlin modules.

 This construction does not work if  $\Im\Phi$  is not holomorphic in any neighborhood of $\gimel(v)$. The study of the case when $\Im\Phi$  is not holomorphic in $\gimel(v)$ but has at most one simple pole, or in other words $\Im\Phi$  is $1$-singular, was initiated by V.~Futorny, D.~Grantcharov and E.~Ramirez  in \cite{Futorny}. The authors \cite{Futorny} constructed the universal $1$-singular Gelfand-Tsetlin $\mathfrak{gl}_n(\mathbb C)$-module using additional formal variables that were called {\it derivative tableaux}. For another construction of the universal $1$-singular Gelfand-Tsetlin $\mathfrak{gl}_n(\mathbb C)$-module see \cite{Zad}, which was posted to the arXiv when the present paper was in preparation.

In the present paper we define a subring $\mathcal D_v$ of $\mathcal R$, where $v$ is a certain point of a $1$-singular $\gimel$-orbit. To the ring $\mathcal D_v$ we associate  the vector space $\mathcal S=\mathcal S(\mathcal D_v)$ with  basis $\mathcal B = \mathcal B(\mathcal D_v)$ consisting of some local distributions supported at $\gimel(v)$ such that $\mathcal S$ is a natural $\mathcal D_v$-module. In particular this implies the following universal property of $ \mathcal D_v$: for any homomorphism of rings $\Psi:	\mathcal U(\mathfrak h) \to \mathcal D_v$ the vector space $\mathcal S$ is also an $\mathfrak h$-module. Due to this we call the ring $ \mathcal D_v$ the {\it universal ring}. Further, we observe that $\Phi	(\mathcal U(\mathfrak g)) \subset \mathcal D_v$. Hence our construction gives rise to a $\mathfrak g$-module structure on $\mathcal S$ that is isomorphic to the universal $1$-singular Gelfand-Tsetlin module obtained in \cite{Futorny}. Our observation leads to a new geometric interpretation of the universal $1$-singular Gelfand-Tsetlin module from \cite{Futorny} that allows to simplify proofs from \cite{Futorny} and avoid the use of formal variables.  
 Moreover, similar ideas that we present here can be used in the case of other singularities, see \cite{EMV}.

\bigskip

\textbf{Acknowledgements:} E.~V. was partially partially supported by SFB TR 191 and by the Universidade Federal de Minas Gerais.

\section{Preliminaries}

A Gelfand-Tsetlin tableau is a tableau $(a_{ki})$ of complex numbers, where $1\leq i \leq k \leq n$ and $n\geq 2$. Further we will consider the set $V$ of all Gelfand-Tsetlin tableaux as a complex manifold that is isomorphic to $\mathbb C^{n(n+1)/2}$. 
 Let $\gimel \simeq \mathbb Z^{n(n-1)/2}$ be the free abelian group generated by $\sigma_{st}$, where $1\leq t \leq s \leq n-1$.  We fix the following action of $\gimel$ on $V$: $\sigma_{st} (x) = (x_{ki} +\delta_{ki}^{st})$, where $x = (x_{ki})\in V$ and $\delta_{ki}^{st}$ is the  Kronecker delta. This is $\delta_{ki}^{st}=1$ if $(ki) = (st)$ and $\delta_{ki}^{st}=0$ otherwise.   Further we put $G=S_1\times S_2\times \cdots \times S_n$, so $G$ is the product of symmetric groups $S_i$. The group $G$ acts on $V$ in the following way $(s (x))_{ki} = x_{k s_k(i)},$
where $s=(s_1,\ldots, s_n)\in G$. Denote by $\mathcal M$ and by $\mathcal O$ the sheaves of meromorphic and holomorphic functions on $V$, respectively. Let us take $f\in H^0(V,\mathcal M)$, $s\in G$ and $\sigma\in \gimel$. We set 
$$
 \sigma(f) = f\circ \sigma^{-1}, \quad s(f) = f\circ s^{-1},\quad s(\sigma) = s\circ \sigma\circ s^{-1}.
$$ 
These formulas define an action of $\gimel$ on $H^0(V,\mathcal M)$ and actions of $G$ on $H^0(V,\mathcal M)$ and on $\gimel$, respectively.

Denote by $\mathcal M\star \gimel$ the sheaf of meromorphic functions on $V$ with values in $\gimel$. An element of $H^0(V,\mathcal M\star \gimel)$ is a finite sum $\sum\limits_i f_i \sigma_i$, where $f_i\in \mathcal M$ and $\sigma_i\in \gimel$. In other words, $\mathcal M\star \gimel$ is the sheaf of meromorphic sections of the trivial bundle $V\times \bigoplus_{\sigma\in \gimel} \mathbb C \sigma \to V$. There exists a structure of a skew group ring on $H^0(V,\mathcal M\star \gimel)$, see \cite{FO}. Indeed, 
$$
\sum_if_i \sigma_i \circ 
\sum_j f'_j \sigma'_j := 
\sum_{ij} f_i \sigma_i(f'_j) \sigma_i\circ \sigma'_j.
$$
Here $f_i, f'_j\in H^0(V,\mathcal M)$  and $\sigma_i,\sigma'_j \in \gimel$. This skew group ring we denote by $\mathcal R$. To simplify notations we use $\circ$ for the multiplication in $\mathcal R$ and for the product in $\gimel$. 
 On $H^0(V,\mathcal M\star \gimel)$ we will consider also  the following multiplication $A*B := B\circ A$.

Recall that a Gelfand-Tsetlin tableau is called {\it generic} if $x_{rt} - x_{rs} \notin \mathbb Z$ for any $r$ and  $s\ne t$. The definition of a standard Gelfand-Tsetlin tableau can be found in \cite{Futorny}.
 The classical Gelfand-Tsetlin formulas  have the following form in terms of generators of $\mathfrak{gl}_n(\mathbb C)$, see for instance \cite{Futorny}, Theorems $3.6$ and $3.8$, for details. 
\begin{equation}\label{eq G-Ts generators}
\begin{split}
	E_{k,k+1}(v) &= - \sum_{i=1}^k \frac{\prod_{j= 1}^{k+1} (x_{ki}- x_{k+1,j})}{\prod_{j\ne i}^k (x_{ki}- x_{kj})} (v + \delta_{ki});\\
	E_{k+1,k} (v) & = \sum_{i=1}^k \frac{\prod_{j= 1}^{k-1} (x_{ki}- x_{k-1,j})}{\prod_{j\ne i}^k (x_{ki}- x_{kj})} (v - \delta_{ki});\\
	E_{k,k} (v) & = \Big( \sum_{i=1}^k (x_{ki} +i-1) - \sum_{i=1}^{k-1} (x_{k-1,i} +i-1)  \Big) (v),
\end{split}
\end{equation}
 Here $E_{st}\in \mathfrak{gl}_n(\mathbb C)$ are standard generators and $v\in V$ is either a standard or generic Gelfand-Tsetlin tableau with coordinates $v=(x_{ki})$ and $(v \pm \delta_{ki}) = \sigma^{\pm 1}_{ki}(v)$.  Assume that $v$ is a generic Gelfand-Tsetlin tableau. Theorem $3.8$ in \cite{Futorny} says that Formulas (\ref{eq G-Ts generators}) define a $\mathfrak g$-module structure on the vector space spanned by the elements of the orbit $\gimel(v)$.

Let us identify the point $v\in V$ with the evaluation map $ev_v: H^0(V,\mathcal O)\to \mathbb C$. Then we can  define the map  $E_{st} \mapsto \Phi(E_{st})\in \mathcal R$ using the equality $ev_{v} \circ \Phi(E_{st})  = E_{st} (v)$ for $v$ generic. Since generic points are dense in $V$,  $\Phi(E_{st})$ is a well-defined element of $\mathcal R$. For example,
 $$
 E_{k,k+1} = - \sum_{i=1}^k \frac{\prod_{j= 1}^{k+1} (x_{ki}- x_{k+1,j})}{\prod_{j\ne i}^k (x_{ki}- x_{kj})} \sigma^{-1}_{ki}.
 $$
  In \cite{FO} the following theorem was proved.

 \medskip
 
 \t\label{teor Fut Ovs} {\sl The map 
 	$\Phi:	\mathcal U(\mathfrak g) \to \mathcal R$ is a homomorphism of rings. Here $E_{st} \mapsto \Phi(E_{st})$, where  $E_{st}\in \mathfrak{gl}_n(\mathbb C)$, is as above. 
 	
 }
 
 \medskip
 
\noindent{\bf Remark.} Note that $\Im(\Phi)$ is $G$-invariant. This fact can be verified directly. 

 \medskip

From Theorem \ref{teor Fut Ovs} it follows that for any generic $x\in V$ the formula
$ \Phi(E_{st}) (ev_y) = ev_y \circ \Phi(E_{st})$
defines an action of $\mathfrak{gl}_n(\mathbb C)$ on the vector space spanned by local distributions $ev_y$, where $y\in \gimel(x)$. Here $ev_y\circ (f\sigma) = ev_y(f) ev_y\circ \sigma$ and $ev_y\circ \sigma(g) = ev_y(\sigma(g))$ for $g\in \mathcal O$. Since elements $\Phi(E_{st})$ are holomorphic in a sufficiently small neighborhood of the orbit $\gimel(x)$, the expression $ev_y \circ \Phi(E_{st})$ is well-defined. More generally, any homomorphism of rings $\Psi:	\mathcal U(\mathfrak h) \to \mathcal R$, where $\mathfrak h$ is any Lie algebra, such that the image $\Psi(\mathfrak h)$ is holomorphic in a neighborhood of $\gimel(x)$ defines an action of $\mathfrak{gl}_n(\mathbb C)$ on the vector space spanned by the local distributions $ev_y$, where $y\in \gimel(x)$. The interpretation of  a point $y\in V$ as a local distribution $ev_y$ suggests the possibility to define a $\mathfrak{gl}_n(\mathbb C)$-module structure on other local distributions, i.e. on linear maps $D_y: \mathcal O_y\to \mathbb C$ with $\mathfrak m_y^p\subset \Ker(D_y)$, where $p>0$ and  $\mathfrak m_y$ is the maximal ideal in the local algebra $\mathcal O_y$. This idea we develop in the present paper.

The main problem is that the ring $\mathcal R$ does not act on the vector space of local distributions, because of singularities. 
In the next section we will construct the universal ring $\mathcal D_v\subset \mathcal R^{G_v}$, where $v$ is a certain $1$-singular point in $V$ and $G_v\subset G$ is the stabilizer of $v$. We will show that $\mathcal D_v$ acts on $G_v$-invariant holomorphic functions $H^0(V,\mathcal O^{G_v})$, where the action is given by $(f\circ \sigma) (F) = f F \circ \sigma^{-1}$. This action induces an action $(f\circ \sigma)(D_y)  =   D_y \circ (f\circ \sigma)$ of $(\mathcal D_v,*)$ on $G_v$-invariant holomorphic local distributions $D: H^0(V,\mathcal O^{G_v}) \to \mathbb C$ supported at $\gimel(v)$.  By Theorem \ref{teor Fut Ovs} we have also a structure of $\mathfrak{gl}_n(\mathbb C)$-module on the vector space of these local distributions. Further we will consider local distributions $ev_v\circ A : H^0(V,\mathcal O^{G_v}) \to \mathbb C$, where $A\in \mathcal D_v$. Clearly this vector space is $\mathcal D_v$- and hence $\mathfrak{gl}_n(\mathbb C)$-submodule. The last step is to find a basis $\mathcal B_v$ for the vector space spanned by $\{ev_{v} \circ A,\,\,|\,\, A\in \mathcal D_v\}$. This basis we call the universal basis for the universal ring $\mathcal D_v$.

Our construction implies that for any homomorphism $\Psi:	\mathcal U(\mathfrak h) \to \mathcal D_v$, where $\mathfrak h$ is a Lie algebra, $\mathcal B_v$ is a basis for the corresponding $\mathfrak h$-module. We will see that $\Im (\Phi) \subset \mathcal D_v$ and that $\mathcal B_v$ coincide with the basis constructed in \cite{Futorny}. We develop these ideas in the case of any point $x\in V$ in \cite{EMV}.

\section{Main result}

 A point $x=(x_{kj})\in V$ is called {\it $1$-singular} if there exist $x_{ki}$ and $x_{kj}$, where $i\ne j$, such that $x_{ki} - x_{kj}\in \mathbb Z$ and  $x_{rs} - x_{rt}\notin \mathbb Z$, where $s\ne t$, for each $(r,s,t)\ne (k,i,j)$. Note  that the generators from (\ref{eq G-Ts generators}) have one simple pole at the orbit $\gimel(x)$ for any $1$-singular point $x$. Let us fix an $1$-singular point $x^0=(x^0_{kj})\in V$ such that $x^0_{ki} - x^0_{kj}\in \mathbb Z$.  We put $z_1= x_{ki} - x_{kj}$, $z_2=x_{ki} + x_{kj}$ and we denote  by $z_3,z_4,\ldots$ other coordinates $(x_{st})$ in $V$. So $(z_i)$ are new coordinates in $V$. From now on we fix a point $v=(0,z^0_2,\ldots, z^0_n)\in \gimel(x^0)$  and a sufficiently small neighborhood $W$ of the orbit  $\gimel(v) = \gimel(x^0)$ that is invariant with respect to the group $\gimel$ and with respect to $\tau\in G$, where $\tau$ is defined by $\tau(z_1)=-z_1$ and $\tau(z_i)=z_i$, $i>1$.
 From now on we will consider restrictions of elements of $\mathcal R$ on $W$. We denote by $G_v =\{\id,\tau\}\subset G$ the stabilizer of $v$.

We say that an element $A\in\mathcal R$ is at most $1$-singular at $v$, if $A=\sum_ih_i\sigma_i$, where $h_i$ are holomorphic at $v$ or have the form $h_i=g_i/z_1$, where $g_i$ are holomorphic at $v$. We need the following proposition.

\medskip
\prop\label{prop A_1 ...A_n is 1 sing} { \sl Let $A_j = \sum\limits_i (H^j_{i}/z_1)\sigma_{i}\in \mathcal R^{G_v}$, where $j=1,\ldots,m$ and $H^j_{i}$ are holomorphic in $W$.  Then the product $A_1\circ \cdots \circ A_m$ is at most $1$-singular at $v$.
}

\medskip

\noindent{\it Proof.}  Assume by induction that for $k=m-1$ our statement holds. In other words, assume that $A_1\circ \cdots \circ A_{m-1} = \sum\limits_i(G_{i}/z_1)\sigma_{i}$, where $G_{i}$ are holomorphic at $v$.  We have
\begin{align*}
A_1\circ \cdots \circ A_{m} = \sum\limits_{i,j} \frac{G_{i} \sigma_i(H^m_{j})}{z_1\sigma_i(z_1)} \sigma_{i} \circ \sigma_{j}.  
\end{align*}
Assume that this product is two singular at $v$. Note that $\sigma_i(H^m_{j})$ is holomorphic in $W$. Therefore, $\sigma_{i_0}(z_1) = z_1$ for a certain $i_0$ and, hence,  $\tau(\sigma_{i_0}) = \sigma_{i_0}$. Further, $\tau(\sum\limits_i(G_{i}/z_1)\sigma_{i}) = \sum\limits_i(G_{i}/z_1)\sigma_{i}$, since $A_1\circ \cdots \circ A_{n}$ is $\tau$-invariant. Therefore $\tau(G_{i_0}/z_1 \sigma_{i_0}) = \tau(G_{i_0}/z_1) \sigma_{i_0} = G_{i_0}/z_1 \sigma_{i_0}$.  Hence $\tau(G_{i_0}/z_1)$ is also $\tau$-invariant and $\tau(G_{i_0}) = -G_{i_0}$. Therefore, $G_{i_0} = z_1G'_{i_0}$, where $G'_{i_0}$ is holomorphic at $v$. Therefore, $G_{i_0} \sigma_{i_0}(H^m_{j})/ z_1\sigma_{i_0}(z_1) = G'_{i_0} \sigma_{i_0}(H^m_{j})/ z_1$ is $1$-singular. The proof is complete.$\Box$ 

\medskip

\noindent{\bf Remark.} Elements  $A_j = \sum\limits_i (H^j_{i}/z_1)\sigma_{i}$ as in Proposition \ref{prop A_1 ...A_n is 1 sing} generate a subring $\mathcal D_{v}$ in $\mathcal R^{G_v}$. We call this ring the {\it universal ring of} $v$. By  Proposition \ref{prop A_1 ...A_n is 1 sing} any element in $\mathcal D_{v}$ is at most $1$-singular at $v$. If $A = \sum\limits_i (H_{i}/z_1)\sigma_{i}$ is a generator of $\mathcal D_{v}$ and $F\in H^0(W,\mathcal O^{G_v})$, then $A(F) = F'/z_1$ is at most $1$-singular at $v$, holomorphic in $W\setminus\{v \}$  and $G_v$-invariant. Therefore 
 $\tau(F') = -F'$ and hence $A(F)$ is holomorphic. So we defined an action of $\mathcal D_{v}$ on $H^0(W,\mathcal O^{G_v})$.  
 
\medskip

 We put $g_i:= z_1h_i$, where $\sum\limits_ih_i\sigma_i\in \mathcal R$.  Consider the following set of $G_v$-invariant local distributions defined on  $H^0(W,\mathcal O^{G_v})$:
 \begin{equation}\label{eq diff operator basis}
 \begin{split}
 D^1_{\sigma}:=\frac12 ev_v \circ (\sigma + \tau(\sigma) ),\quad
 D^2_{\sigma'}:= ev_v \circ \frac{(\sigma' - \tau(\sigma') )}{2z_1}, \quad \sigma, \sigma' \in \gimel, \,\,\tau(\sigma')\ne \sigma'.
 \end{split}
 \end{equation}
 Note that $(\sigma' - \tau(\sigma'))/ z_1$ and $\sigma + \tau(\sigma)$ are elements of $\mathcal D_{v}$, hence Formulas (\ref{eq diff operator basis}) are well-defined. 
 Moreover we have the following equalities 
 \begin{equation}\label{eq relations n=1}
 D^1_{\tau(\sigma)} = D^1_{\sigma} \quad \text{and} \quad D^2_{\tau(\sigma')} = - D^2_{\sigma'}.
 \end{equation}
Denote $\Delta: = \{\sigma\in \gimel\,\,|\,\, \sigma(x_{ki},x_{ki}) = (x_{ki}+m_1,x_{ki}+m_2), \,\, m_1 \leq m_2 \}$ and consider the set $\mathcal B_v:= \{D^1_{\sigma}, D^2_{\sigma'}\,\, | \,\, \sigma,\sigma'\in \Delta,\,\, \tau(\sigma')\ne \sigma' \}$. The set $\mathcal B_v$ is a set of linearly independent distributions defined on $H^0(W,\mathcal O^{G_v})$. To see this we should apply $D^i_{\sigma}$ to a linear combination $\alpha + \beta z_1^2$ of $G_v$-invariant functions $1$ and $z_1^2$. Hence $\mathcal B_v$ is a basis of the vector subspace $S$ in  $H^0(W,\mathcal O^{G_v})^*$ spanned by elements from $\mathcal B_v$. In the next proposition we show that $S$ is a $\mathcal D_v$- and $\mathfrak {gl}_n(\mathbb C)$-module. This $\mathfrak {gl}_n(\mathbb C)$-module is isomorphic to the universal $1$-singular Gelfand-Tsetlin module constructed  in \cite{Futorny}, see Section $4$ for details.

\medskip

\prop\label{prop formulas n=1} {\sl Let us take $\sum\limits_ih_i\sigma_i\in \mathcal D_v$. Then   we have
\begin{equation}\label{eq formula n=1, sing}
\begin{split}
ev_v \circ (\sum\limits_ih_i\sigma_i) = \sum\limits_i g_i(v) \cdot D^2_{\sigma_i} +
\sum\limits_i\frac{\partial g_i}{\partial z_1}(v)   
D^1_{\sigma_i},
\end{split}
\end{equation}
where $g_i = z_1 h_i$.
Note that in the case if $h_i$ is holomorphic, we have $g_i(v)=0$ and $\frac{\partial g_i}{\partial z_1}(v) = h_i(v)$. Therefore, $S$  is a $\mathcal D_v$-module with basis $\mathcal B_v$. 

}

\medskip
\noindent{\it Proof.} Using the series expansion $g_i = g_i|_{z_1=0} + \frac{\partial g_i}{\partial z_1}|_{z_1=0} z_1 + \cdots$, we get
$$
ev_v \circ \sum\limits_ih_i\sigma_i = \sum\limits_i g_i(v) ev_v \circ \frac{\sigma_i}{z_1} + \sum\limits_i \frac{\partial g_i}{\partial z_1}(v) ev_v \circ \sigma_i.
$$
Note that $ev_v \circ z_1^m\sigma_i =0$ for $m>1$. Using the symmetrization $2ev_v \circ \sum\limits_ih_i\sigma_i = ev_v \circ \sum\limits_ih_i\sigma_i + \tau(ev_v \circ \sum\limits_ih_i\sigma_i)$, we obtain the result.$\Box$ 

\medskip

Let $X$ be one of generators (\ref{eq G-Ts generators}). Clearly $\Phi(X)\in \mathcal D_v$, see Remark after Theorem \ref{teor Fut Ovs}.

\medskip

\t\label{teor main}{\bf [Main result 1]} {\sl 
	The vector space $S$ spanned by elements of $\mathcal B$ is a $\mathfrak{gl}_n(\mathbb C)$-module. The action is given by Formulas (\ref{eq formula n=1, sing}).
	
}

\medskip
\noindent{\it Proof.}  The result follows from Theorem \ref{teor Fut Ovs} and Proposition \ref{prop formulas n=1}. Indeed, $\Im(\Phi)\subset \mathcal D_v$, hence we get a $\mathfrak{gl}_n(\mathbb C)$-module.$\Box$

\medskip

 In fact we proved a more general result than it is formulated in  Theorem \ref{teor main}.

\medskip

\t\label{teor main general}{\bf [Main result 2]} {\sl Let $\mathfrak h$ be a Lie algebra and $\Psi : \mathcal U(\mathfrak h) \to \mathcal D_v$ be a homomorphism of rings. Then the vector space $S$ spanned by elements from  $\mathcal B$  is an $\mathfrak h$-module. In other words the basis $\mathcal B$ is universal for any homomorphisms $\Psi : \mathcal U(\mathfrak h) \to \mathcal D_v$. 
	
}

\section{Appendix} 

Theorem \ref{teor main} recovers one of the main results of \cite{Futorny}, a construction of the universal $\mathfrak{gl}_n(\mathbb C)$-module. Another main result of \cite{Futorny} is that in many cases  the $\mathfrak{gl}_n(\mathbb C)$-module $S$ is irreducible, see Theorem $4.14$. Let us give an explicit correspondence between the basis constructed in \cite{Futorny} and our basis $\mathcal B$. We use notations from \cite{Futorny}. In \cite{Futorny} the authors consider the basis $\{ T(\sigma(v)),\,\, \mathcal {D}T(\sigma'(v)) \}$, where $\sigma,\sigma'\in \gimel$, such that $T(\sigma(v)) - T(\tau(\sigma(v)))=0$, $ \mathcal {D}T(\sigma'(v)) +  \mathcal {D}T(\tau(\sigma'(v)))=0$ and $\tau(\sigma')\ne \sigma'$, see Remark $4.5$ in \cite{Futorny}. The element $T(\sigma(v))$ was considered as a point in $V$ and  $\mathcal {D}T(\sigma'(v))$ was considered as a formal additional variable. (In our notations, $T(\sigma(v))$ is just $\sigma(v)\in V$, where $v$ as above.) Further, in \cite{Futorny} the action of $\mathfrak{gl}_n(\mathbb C)$ is given by the following formulas,  \cite[Theorem 4.11]{Futorny}: 
\begin{align*}
E_{rs} (T(\sigma(v)))& = ev_{v}\circ\frac{\partial}{\partial z_1}  z_1 E_{rs}(T(\sigma(x)));\\
E_{rs} (\mathcal {D}T(\sigma'(v))& = ev_{v}\circ\frac{\partial}{\partial z_1}   E_{rs}(T(\sigma'(x))), \,\,\, E_{rs}\in \mathfrak{gl}_n(\mathbb C),
\end{align*}
where $x=(x_{ki})$ are coordinates in a neighborhood of $v$.  
The explicit correspondence between the bases is given by the following formulas   
\begin{align*}
&2D^2_{\sigma}(T(v)) =  \mathcal {D}T(\sigma'(v)) - \mathcal {D}T(\tau(\sigma'(v))), \quad
 2D^1_{\sigma}(T(v)) =  T(\sigma(v)) + T(\tau(\sigma(v))).
\end{align*}

\bigskip

\noindent
E.~V.: Departamento de Matem{\'a}tica, Instituto de Ci{\^e}ncias Exatas,
Universidade Federal de Minas Gerais,
Av. Ant{\^o}nio Carlos, 6627, Caixa Postal: 702, CEP: 31270-901, Belo Horizonte, 
Minas Gerais, BRAZIL,
email: {\tt VishnyakovaE\symbol{64}googlemail.com}

\end{document}